\documentclass[11pt]{article}
\usepackage{epic,latexsym,amssymb}
\usepackage{amsfonts}
\usepackage{amscd}
\usepackage{amsmath}
\usepackage{graphicx}
\usepackage{color}
\usepackage{caption,
caption}
\usepackage{tikz}
\usepackage{amsthm}
\usepackage{bbm}
\usepackage{amsfonts}

\usepackage[margin=0.25in]{geometry}
\usepackage{pgfplots}
\pgfplotsset{width=10cm,compat=1.9}

\usepgfplotslibrary{external}
\usepackage{epic,latexsym,amssymb}
\usepackage{tikz}
\usepackage[bold,full]{complexity}
\usepackage{tikz,epsf}
\usepackage{hyperref}
\usepackage{graphicx}[width=\textwidth]
\usepackage{amsmath}
\usepackage{listings}

\usepackage{amsfonts}
\usepackage{amscd}
\usepackage{amsmath}
\usepackage{graphicx}
\usepackage{color}
\usepackage{caption,subcaption}
\usepackage{lineno}

\usepackage[ruled,vlined]{algorithm2e}

\usepackage{lineno}
\newtheorem{thma}{Theorem}

\newtheorem{thm}{Theorem}
\newtheorem{conj}[thm]{Conjecture}
\newtheorem{ob}[thm]{Observation}

\newtheorem{lem}[thm]{Lemma}
\newtheorem{quest}{Question}

\newtheorem{cor}[thm]{Corollary}

\newcommand{\1}{\vspace{0.1cm}}

\newcommand{\QEDmark}{\mbox{\textsc{qed}}}
\newcommand{\proofStarter}[1]{\textsc{#1} }

\def\vertex(#1){\put(#1){\circle*{2}}}
\def\vertexo(#1){\put(#1){\circle{2}}}
\def\vert(#1){\put(#1){\circle*{1.5}}}
\def\verto(#1){\put(#1){\circle{1.5}}}
\def\lab(#1)#2{\put(#1){\makebox(0,0)[c]{#2}}}
\definecolor{DarkGreen}{rgb}{0.2, 0.6, 0.3}

\definecolor{electricindigo}{rgb}{0.44, 0.0, 1.0}

\let\oldenumerate\enumerate
\renewcommand{\enumerate}{
  \oldenumerate
  \setlength{\itemsep}{0.5pt}
  \setlength{\parskip}{0pt}
  \setlength{\parsep}{0pt}
}

\begin{document}

\title{Computer assisted discovery: Zero forcing vs vertex cover}
\author{$^{1}$Boris Brimkov, $^{2}$Randy Davila\footnote{Corresponding author}, $^{3}$Houston Schuerger, and $^{4}$Michael Young\\
\\
$^1$Department of Mathematics and Statistics\\
Slippery Rock University\\
Slippery Rock, PA 16057, USA\\
\small {\tt Email: boris.brimkov@sru.edu} \\
\\
$^2$Research and Development \\
RelationalAI \\
Berkeley, CA 94704, USA\\
\small {\tt Email: randy.davila@relational.ai}\\
\\
$^3$Department of Mathematics\\
Trinity College\\
Hartford, CT 06106, USA\\
\small {\tt Email: houston.schuerger@trincoll.edu} \\
\\
$^4$Department of Mathematical Sciences\\
Carnegie Mellon University\\
Pittsburgh, PA 15213, USA \\
\small {\tt Email: michaely@andrew.cmu.edu} \\
\\
}

\date{}
\maketitle

\begin{abstract}
In this paper, we showcase the process of using an automated conjecturing program called \emph{TxGraffiti} written and maintained by the second author. We begin by proving a conjecture formulated by \emph{TxGraffiti} that for a claw-free graph $G$, the vertex cover number $\beta(G)$ is greater than or equal to the zero forcing number $Z(G)$. Our proof of this result is constructive, and yields a polynomial time algorithm to find a zero forcing set with cardinality $\beta(G)$. We also use the output of \emph{TxGraffiti} to construct several infinite families of claw-free graphs for which $Z(G)=\beta(G)$. Additionally, inspired by the aforementioned conjecture of \emph{TxGraffiti}, we also prove a more general relation between the zero forcing number and the vertex cover number for any connected graph with maximum degree $\Delta \ge 3$, namely that $Z(G)\leq (\Delta-2)\beta(G)$+1. 

%\Randy{
%\emph{TxGraffiti} is an open source automated conjecturing program written in Python (versions 3.6 and higher), written and maintained by the second author. In recent years this program has produced a myriad of conjectures in graph theory leading to research publications. Most notably, TxGraffiti has produced surprising conjectures on graph invariants such as the \emph{independence number}, the \emph{domination number}, the \emph{matching number}, and the \emph{zero forcing number}. The conjectures of TxGraffiti are posed in the form of inequalities between graph invariants, and in this paper we prove one of these  conjectures which relates the zero forcing number and the vertex cover number of a simple and connected graph. Specifically, we prove that if $G$ is a connected and claw-free graph, then $Z(G) \le \beta(G)$, where $Z(G)$ and $\beta(G)$ denote the zero forcing number and vertex cover number of $G$, respectively. This bound is sharp and improves on known bounds for the zero forcing number. Additionally, and inspired by the aforementioned conjecture of TxGraffiti, we also prove a general theorem bounding the zero forcing number by the vertex cover number. Specifically, if $G$ is a connected graph, and $\mathbbm{1}_{bipartite}(G)$ denotes the indicator variable for $G$ being a bipartite graph, and $\Delta \ge 3$ is the maximum degree of $G$, then $Z(G) \le (\Delta - 2)\beta(G) + \mathbbm{1}_{bipartite}(G)$. This inequality is also sharp, and we provide examples achieving this. }
\end{abstract}


\section{Introduction}

Data-driven machine learning techniques have become ubiquitous in various scientific fields and often yield unexpected results such as DeepMind's protein folding solution~\cite{protein} or AlphaGo's novel game strategies~\cite{AlphaGo}. Such techniques can also be used in pure mathematics to come up with new and exciting computer generated conjectures. The idea of automated conjecturing was first discussed by Turing~\cite{Turing} in the 1950s. Over the next three decades, there were several attempts to implement automated conjecturing programs, but the resulting attempts generated thousands of conjectures and were difficult to parse for meaningful problems. The first program to circumvent this problem and produce conjectures suitable for publication was Fajtlowicz's GRAFFITI~\cite{Graffiti}, named so because its conjectures were ``written on the wall'' for other mathematicians to view. The conjectures of GRAFFITI attracted the attention of many well-known mathematicians including Erd\"{o}s, Chung, and Faber. Subsequently, several other viable automated conjecturing programs were developed, such as DeLaVi\~{n}a's GRAFFITI.pc~\cite{Graffiti.pc}, Larson's Conjecturing~\cite{Larson}, Lenat’s AM~\cite{Lenat_1, Lenat_2, Lenat_3}, Epstein’s GT~\cite{Epstein_1, Epstein_2}, Colton’s HR~\cite{Colton_1, Colton_2, Colton_3}, Hansen and Caporossi’s AGX~\cite{AGX_1, AGX_2, AGX_3}, Mélot’s Graphedron~\cite{graphedron_1}, and Davila's \emph{TxGraffiti}~\cite{TxGraffiti} and \emph{Conjecturing.jl}~\cite{Conjecturing.jl}. These programs generate conjectures about various mathematical objects, including matrices, numbers, graphs, and functions.

\begin{figure}[htb]
\begin{center}
\includegraphics{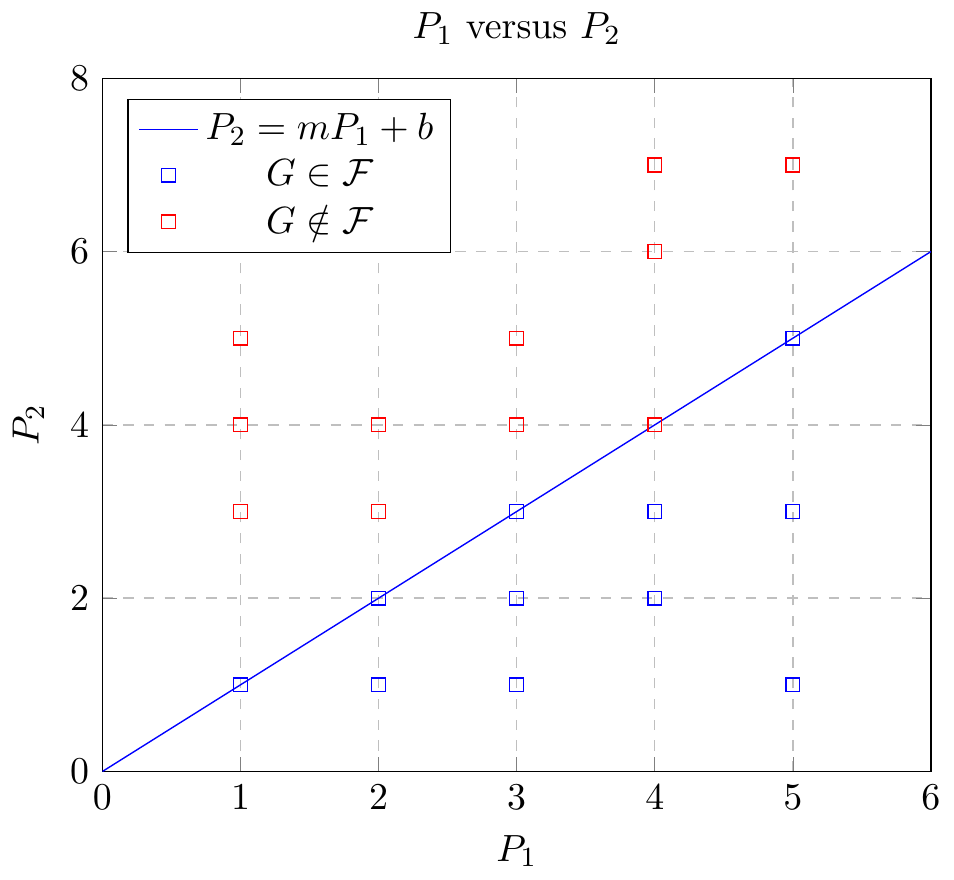}
\end{center}
\caption{Illustration of the TxGraffiti automated conjecturing process. }
\label{fig:construction0}
\end{figure}

In this paper, we showcase the utility of automated conjecturing by proving a conjecture generated by \emph{TxGraffiti}, and using that proof as a springboard for a more general result. Inspired by GRAFFITI and GRAFFITI.pc, \emph{TxGraffiti} was written in 2017 by the second author using \emph{Python} (versions 1.6 and higher). This program generates conjectures in the form of inequalities relating graph parameters and has already led to several new results and publications (see, e.g.,~\cite{CaDaHePe2022, DaHe19b, DaHe21a, DaHe19c}). Broadly speaking, the program computes various graph parameters for different families of graphs that are stored in a large database and outputs relationships that are not violated by any of the graphs in a given family. A simple illustration of the design of \emph{TxGraffiti} is shown in Figure~\ref{fig:construction0}. Two graph parameters, $P_1(G)$ and $P_2(G)$, are computed for a number of graphs, some of which belong to a family $\mathcal{F}$. The program then solves a linear program and finds the variables $m$ and $b$ that satisfy the line $P_1 = mP_2 + b$ which separates the graphs according to membership in $\mathcal{F}$. Then, it outputs the conjecture ``If $G\in \mathcal{F}$, then $P_1(G)\leq mP_2(G) + b$." The conjectures of TxGraffiti are ranked based on a combination of factors, such as whether a proposed inequality subsumes another inequality for a smaller family of graphs. Moreover, a conjecture is only presented if the proposed inequality is sharp for a significant portion of graph instances satisfying the hypotheses of the conjecture.

A particularly strong conjecture made by TxGraffiti was that the vertex cover number of a claw-free graph is an upper bound on the zero forcing number of the graph (see Section \ref{preliminaries} for exact definitions of these parameters). As the first main result of this paper, we prove this conjecture. 
\begin{conj}[TxGraffiti 2019]
\label{conj:main}
If $G$ is a nontrivial connected claw-free graph, then $Z(G) \le \beta(G)$.
\end{conj}
Conjecture \ref{conj:main} is appealing for several reasons. First, it gives a simple relationship between two widely studied graph parameters in an important class of graphs. Moreover, by analyzing the output of \emph{TxGraffiti} in generating Conjecture \ref{conj:main}, we constructed infinite families of claw-free graphs satisfying $Z(G)=\beta(G)$.  Furthermore, our proof of Conjecture~\ref{conj:main} is constructive and yields a polynomial time algorithm to find a zero forcing set of $G$ with cardinality $\beta(G)$. This is useful because the vertex cover number of a claw-free graph can be computed in polynomial time~\cite{Stable-NP}, while the complexity of computing the zero forcing number of a claw-free graph is unknown (and is NP-hard in general~\cite{NP-Complete}). Finally, inspired by the proof technique for Conjecture \ref{conj:main}, we also discovered the following more general relation between the zero forcing number and the vertex cover number of an arbitrary graph. 
\begin{thm}
\label{theorem2a}
If $G$ is a connected graph with maximum degree $\Delta \geq 3$, then 
$$
Z(G) \le (\Delta - 2)\beta(G) + 1,
$$
and this bound is sharp.
\end{thm}

This paper is organized as follows. In Section \ref{preliminaries} we recall some graph theoretic terminology and notation. In Section \ref{sec:thm1} we give the proof of Conjecture~\ref{conj:main}; we also provide an approximation algorithm for zero forcing in claw-free graphs and give constructions of graphs with zero forcing number equal to the vertex cover number. In Section~\ref{sec:thm2} we give the proof of Theorem~\ref{theorem2a} and provide families of graphs for which the bound in Theorem~\ref{theorem2a} holds with equality. We conclude with some final remarks and open questions in Section~\ref{conclusion}.

\section{Preliminaries}\label{preliminaries}
Throughout this paper, all graphs considered will be simple, undirected, and finite. Let $G$ be a graph with vertex set $V(G)$ and edge set $E(G)$. The \emph{order} of $G$ is $n(G) = |V(G)|$. Two vertices $v,w \in V(G)$ are \emph{neighbors}, or \emph{adjacent}, if $vw \in E(G)$. The \emph{open neighborhood} of $v\in V(G)$, is the set of neighbors of $v$, denoted $N_G(v)$; the \emph{closed neighborhood} of $v$ is $N_G[v] = N_G(v) \cup \{v\}$. The closed neighborhood of $S\subseteq V$ is $N_G[S] = \bigcup_{v\in S}N_G[v]$. The \emph{degree} of a vertex $v\in V(G)$, denoted $d_G(v)$, is equal to $|N_G(v)|$. The maximum and minimum degree of $G$ will be denoted $\Delta(G)$ and $\delta(G)$, respectively. When there is no scope for confusion, we will use the notation $n = n(G)$, $\delta = \delta(G)$, and $\Delta = \Delta(G)$, to denote the order, minimum degree, and maximum degree, respectively. 

Two vertices in a graph $G$ are \emph{independent} if they are not neighbors. A set of pairwise independent vertices in $G$ is an \emph{independent set} of $G$. The number of vertices in a maximum independent set in $G$ is the \emph{independence number} of $G$, denoted $\alpha(G)$. A \emph{vertex cover} for $G$ is a set of vertices $C\subseteq V(G)$ so that every edge in $G$ has one endpoint in $C$. The \emph{vertex cover number} of $G$, denoted $\beta(G)$, is the minimum cardinality of a vertex cover in $G$. Note that for any graph $G$, $\beta(G) = n(G) - \alpha(G)$. 

We denote the cycle and complete graph on $n$ vertices by $C_n$ and $K_n$, respectively. A \emph{leaf} is a vertex of degree one, while its neighbor is a \emph{support vertex}. A \emph{strong support vertex} is a vertex with at least two leaf neighbors. A star is a non-trivial tree with at most one vertex that is not a leaf. Thus, a \emph{star graph} is the tree $K_{1,k}$ for some $k \ge 1$. A graph $G$ is \emph{$F$-free} if $G$ does not contain $F$ as an induced subgraph. In particular, if $G$ is $F$-free, where $F = K_{1,3}$, then $G$ is \emph{claw-free}. Claw-free graphs have been widely studied and a comprehensive survey of claw-free graphs has been written by Flandrin, Faudree, and Ryj{\'a}{\v{c}}ek ~\cite{claw-free}.

The \emph{zero forcing process} on $G$ is defined as follows: Let $B\subseteq V(G)$ be a set of initially ``blue colored'' vertices, all remaining vertices being ``white colored''. At each discrete time step, if a blue colored vertex has a unique white colored neighbor, then this blue colored vertex \emph{forces} its white colored neighbor to become colored blue. If $v$ is a blue vertex which forces a white colored neighbor to be colored blue, then we say that $v$ has been \emph{played}. The initial set of blue colored vertices $B$ is a \emph{zero forcing set}, if by iteratively applying the zero forcing process all of $V(G)$ becomes colored blue. The \emph{zero forcing number} of $G$, written $Z(G)$, is the cardinality of a minimum zero forcing set in $G$. Zero forcing was introduced in \cite{AIM-Workshop} as a bound on the minimum rank over all symmetric matrices that have the same off-diagonal nonzero pattern as the adjacency matrix of a given graph. Zero forcing is also related to other processes that arise from the fact that knowing all-but-one of the variables in a linear equation implies the value of the last remaining variable. In particular, processes that are equivalent or similar to zero forcing were independently introduced in quantum control theory \cite{quantum1}, graph searching \cite{fast_mixed_search}, and PMU placement \cite{powerdom3}. See~\cite{Davila2,DaHe18b,DaKaSt18,FuRa19, Genter1, Genter2,FeGrKa19,FeKaSt19, LuTang} for some recent structural results and bounds on the zero forcing number.

For notation and graph terminology not mentioned here we refer the reader to West~\cite{West}.
  %We also will make use of the standard notation $[k] = \{1,\ldots,k\}$.
%%%%%%%%%%%%%%%%%%%%%%%%%%%%%%%%%%%%%%%%%%

\section{Proof of Conjecture \ref{conj:main}}\label{sec:thm1}
In this section we prove Conjecture \ref{conj:main}, which is restated with the following theorem.

\begin{thma}
\label{thm1}
If $G$ is a connected and claw-free graph with $\delta(G) \ge 1$, then 
\[
Z(G) \le \beta(G),
\]
and this bound is sharp.  
\end{thma}

\proof Let $G$ be a connected claw-free graph with minimum degree $\delta(G) \ge 1$. Suppose $C\subseteq V(G)$ is a minimum vertex cover of $G$, and so, $X = V(G)\setminus C$ is a maximum independent set of $G$. Thus, each vertex in $C$ has at least one neighbor in $X$. Furthermore, each vertex in $C$ has at most two neighbors in $X$, since otherwise there would be a vertex in $C$ which together with three of its neighbors in $X$ would induce a claw in $G$. Moreover, since $\delta(G)\ge 1$, each vertex in $X$ is adjacent with at least one vertex in $C$. With these observations, we note that $N_G[C] = V(G)$. Next let $C_1$ denote the set of vertices in $C$ with exactly one neighbor in $X$ and let $C_2$ denote the set of vertices in $C$ with exactly two neighbors in $X$. With the following we construct a zero forcing set of $G$, denoted by $S$, such that $|S| = |C|$ and $|S \cap C| \geq |C|-1$. 

\noindent\textbf{Initialization Phase.} If $C_1 \neq \emptyset$, then we choose $S = C$ as our initial set of blue colored vertices, all other vertices being colored white, which implies all vertices in $X$ are initially colored white. Since $C_1 \neq \emptyset$, observe that if we start the zero forcing process at $S$, then there is at least one vertex in $S$ with exactly one white neighbor, namely some vertex in $C_1$ that is colored blue and has one white neighbor which is contained in $X$. Let this vertex in $C_1$ force its one white neighbor to become colored blue. Hence, each vertex in $C$ is colored blue and at least one vertex in $C$ will have its closed neighborhood colored blue. Furthermore, $|S| = |C|$ and $|S \cap C| \geq |C|-1$. 

If $C_1 = \emptyset$, then there is at least one vertex $v\in C_2$ with two neighbors in $X$, say $u$ and $w$. In this case, we choose $S = (C\setminus\{v\})\cup\{u\}$ as our initial set of blue colored vertices, all other vertices being colored white, which implies all vertices in $X\setminus\{u\}$ are initially colored white. Next observe that if we start the zero forcing process at $S$, then $v$ is the only white colored neighbor of $u$ (which is blue). Thus, $u$ may force $v$ to become colored blue. After $v$ becomes colored blue, then $w$ would be the only white colored neighbor of $v$, and so, $v$ may then force $w$ to become colored blue. Hence, each vertex in $C$ is colored blue and at least one vertex in $C$ will have its closed neighborhood colored blue. Furthermore, $|S| = |C|$ and $|S \cap C| \geq |C|-1$.

\noindent\textbf{Zero Forcing Phase.} Starting from $S$ as our initial set of blue colored vertices, all other vertices being colored white, we allow the zero forcing process to start and continue until no further color changes are possible. Note that in both possible choices of $S$ described above, we are ensured that each vertex in $C$ is colored blue and at least one vertex in $C$ has its closed neighborhood colored blue. Furthermore, note that if all vertices in $N_G[C]$ are colored blue then all vertices in $V(G)$ are colored blue. 

Let $C'$ be the set of vertices in $C$ whose closed neighborhoods are not colored blue after the zero forcing process terminates. Because no further color changes are possible, and because each vertex in $C$ is blue, all vertices in $C'$ have exactly two white neighbors in $X$. This observation implies that if $v\in C'$, then $v\in C_2$. Furthermore, if $v\in C'$, then no vertex of $C\setminus C'$ can be adjacent with a neighbor of $v$ in $X$, since then $v$ would have one or less white neighbors. Moreover, if $v \in C'$, then $v$ cannot be adjacent with $z \in C\setminus C'$, since $G$ being claw-free implies $z$ would have shared at least one neighbor in $X$ with $v$, which is impossible since $z$ would have forced one of the neighbors of $v$ to become colored blue. It follows that no vertex in $N_G[C']$ can be adjacent with a vertex in $C\setminus C'$. Since $G$ is connected and since $C\setminus C' \neq \emptyset$, it must be that $C' = \emptyset$. Thus, $N[C]$ is necessarily colored blue and $S$ was a zero forcing set of $G$. Hence, $Z(G) \le |S| = |C| = \beta(G)$. 

To see this bound is sharp, see Section~\ref{sec:discuss}.  \qed 

The proof of Theorem \ref{thm1} implies that if $G$ is claw-free, then either a minimum vertex cover is a zero forcing set of $G$, or a slight modification of a minimum vertex cover is a zero forcing set of $G$. The precise instructions are given with the following algorithm. Note that finding a maximum independent set in a claw-free graph can be done in polynomial time~\cite{Stable-NP}.

\begin{algorithm}[H]
\SetAlgoLined
\textbf{Input:} Connected claw-free graph $G$ with order $n$ and independence number $\alpha$.
\\
\KwResult{A zero forcing set $B$ of $G$ with cardinality $n - \alpha$.}

Find a maximum independent set $X \subseteq V(G)$\; 
Let $C = V(G)\setminus X$\;
Let $C_1$ be the set of vertices in $C$ with exactly one neighbor in $X$\;
Let $C_2$ be the set of vertices in $C$ with exactly two neighbors in $X$\;
 
 \eIf{$C_1 \neq \emptyset$}{
   $B = C$\;
   }{
   Choose $v\in C_2$ and $w \in N_G(v)\cap X$\;
   $B = (C\setminus\{v\}) \cup \{w\}$\;
  }
 \caption{Zero forcing approximation in claw-free graphs}
\textbf{Return:} $B$\;
\end{algorithm}

\subsection{Sharp Examples of Theorem \ref{thm1}}\label{sec:discuss}
In this section we give infinite families of claw-free graphs with $Z(G) = \beta(G)$. 

\begin{figure}[htb]
\begin{center}
\includegraphics{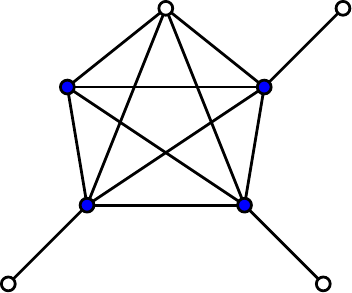}
\caption{One possible claw-free graph $G \in \mathcal{G}'$ with $Z(G) = \beta(G)$. A minimum zero forcing set colored blue; this is also a minimum vertex cover. }
\label{fig:construction1}
\end{center}
\end{figure}

\noindent\textbf{Construction 1.} Let $\mathcal{G}'$ be the family of graphs obtained by starting with the complete graph $K_n$, with $n\ge 3$, and then attach one pendant vertex to at most $n-1$ vertices from $K_n$ to form a graph $G$. For these graphs, we may form a minimum zero forcing set $S$, by letting $S$ be any collection of $n-1$ vertices from $G$ in $K_n$, so that at least one vertex from $S$ is not adjacent with a pendant. This is also a minimum vertex covering for $G$. See Figure~\ref{fig:construction1} for one such construction. 

\begin{figure}[htb]
\begin{center}
\includegraphics{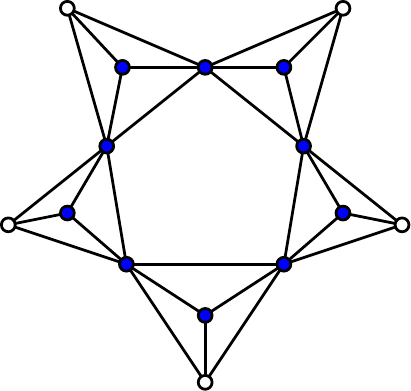}
\end{center}
\caption{One possible claw-free graph $G \in \mathcal{G}''$ with $Z(G) = \beta(G)$. A minimum zero forcing set is colored blue; this is also a minimum vertex cover. }
\label{fig:construction2}
\end{figure}

\noindent\textbf{Construction 2.} Let $\mathcal{G}''$ be the family of graphs obtained by starting with the cycle graph $C_k$, with $k\ge 3$, vertex set $V(C) = \{v_1, v_2, \dots, v_k\}$, and edge set $E(C) =\{v_1v_2, \dots, v_{k-1}v_k, v_kv_1\}$. Next let $K^{1}_{n_1}, K^{2}_{n_2}, \dots, K^{k}_{n_k}$ denote a collection of $k$ complete graphs, where $K^{i}_{n_i}$ denotes the complete graph with order $n_i \ge 1$. For $i=1, \dots k$, attach each vertex of $K^{i}_{n_i}$ to both $v_i$ and $v_{i+1}$ (with $i = i\mod k$). The resulting graph $G$ has independence number $\alpha(G) = k$, which can be seen by taking one vertex from each of the complete graphs $K^{i}_{n_i}$, for $i = 1, \dots , k$. Let $X$ be such a maximum independent set of $G$. Thus, $S = V(G)\setminus X$ is a minimum vertex cover for $G$. Moreover, $S$ is a minimum zero forcing set of $G$. For one example of this construction, see Figure~\ref{fig:construction2}. 

Theorem \ref{thm1} clearly improves on the trivial bound $Z(G) \le n(G) -1$, but can also sometimes improve on other stronger bounds. For example, recall the following upper bound on $Z(G)$ given by Caro and Pepper in~\cite{dynamic-k-forcing}.
\begin{thm}[Caro and Pepper~\cite{dynamic-k-forcing}]
If $G$ is a connected graph with order $n$, maximum degree $\Delta \ge 1$, and minimum degree $\delta$, then
\[
Z(G) \le \frac{(\Delta-2)n - (\Delta - \delta) + 2}{\Delta -1},
\]
and this bound is sharp. 
\end{thm}
Let $G \in \mathcal{G}''$ be the graph obtained with a cycle $C_k$ and $k$ complete graphs $K_2$, as described by Construction 2 in Section~\ref{sec:discuss}. Thus, $n = 3k$, $\Delta = 6$, and $\delta = 3$. It follows that 
\[
Z(G) = \beta(G) = 2k < \frac{12k -1}{5} =  \frac{(\Delta-2)n - (\Delta - \delta) + 2}{\Delta -1}. 
\]
Furthermore, since $F_2(G) \le F_1(G) = Z(G)$, we note that Theorem \ref{thm1} improves on a result given by Amos, Caro, Davila, and Pepper in~\cite{k-forcing}; namely, $F_2(G) \le \beta(G)$ whenever $G$ is claw-free.

\begin{figure}[htb]
\begin{center}
\includegraphics{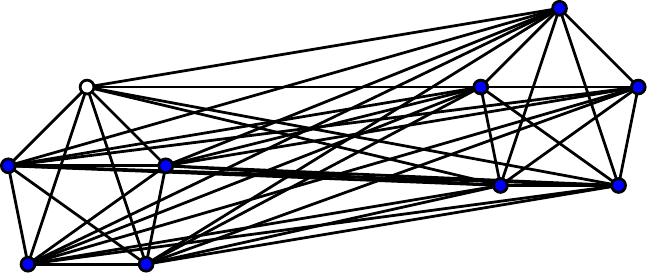}
\end{center}
\caption{The join graph $K_5 \lor K_5$. A minimum zero forcing set colored blue; this is also a minimum vertex cover. }
\label{fig:construction3}
\end{figure}

We next prove a theorem which can be used to iterativey build larger and larger graphs with $Z(G) = \beta(G)$, regardless of whether or not $G$ is claw-free. To do this, recall that the \emph{join graph} of graphs $G$ and $H$, denoted $G\lor H$, is the graph obtained from the disjoint union of $G$ and $H$ by joining each vertex of $G$ to every vertex of $H$. For an illustration see Figure~\ref{fig:construction3}; a minimum zero forcing set is shown in blue. We first recall a useful theorem.

\begin{thm}[Davila, Henning, and Pepper~\cite{DaHePe}]
\label{thmDHP}
If $G$ and $H$ are graphs, then $Z(G\lor H) = \min\{n(G) + Z(H), n(H) +Z(G)\}$.
\end{thm}

With the statement of Theorem \ref{thmDHP} we next provide a construction for graphs with zero forcing number equal to the vertex cover number. Note that this construction can be applied to any graphs $G$ and $H$ that satisfy $Z(G)=\beta(G)$ and $Z(H)=\beta(H)$ regardless of whether $G$ and $H$ are  claw-free.
\begin{thm}\label{thm:join}
If $G$ and $H$ are graphs satisfying $Z(G) = \beta(G)$ and $Z(H) = \beta(H)$, then $Z(G\lor H) = \beta(G\lor H)$. 
\end{thm}
\proof Let $G$ and $H$ be graphs satisfying $Z(G) = \beta(G)$ and $Z(H) = \beta(H)$. Next suppose $X\subseteq V(G) \cup V(H)$ is a maximum independent set in $G\lor H$. Since $v\in V(G)$ and $w\in V(H)$ imply $vw \in E(G\lor H)$, for all $v\in V(G)$ and all $w\in V(H)$, it must be the case that either  $X\subseteq V(G)\setminus V(H)$ or $X\subseteq V(H)\setminus V(G)$. Thus, 
\[
\alpha(G\lor H) = |X| = \max\{\alpha(G), \alpha(H)\},
\]
which implies
\[
\beta(G\lor H) = n(G) + n(H) -  \max\{\alpha(G), \alpha(H)\}.
\]
Next, and without loss of generality, suppose $\alpha(G\lor H) = \alpha(G)$. Since $Z(H) = \beta(H) = n(H) - \alpha(H)$,  Theorem \ref{thmDHP} implies
\begin{equation*}
\begin{split}
Z(G\lor H) & = \min\{n(G)+Z(H), n(H) + Z(G)\} \\
 & = \min\{n(G)+n(H) - \alpha(H), n(H) + n(G) - \alpha(G)\} \\
 & = n(G) + n(H) -  \max\{\alpha(G), \alpha(H)\} \\
 & = \beta(G\lor H),
\end{split}
\end{equation*}
which is the desired equality.\qed 

\section{Proof of Theorem \ref{theorem2}}\label{sec:thm2}
In this section we prove Theorem \ref{theorem2}, which is restated below. 
\begin{thma}
\label{theorem2}
If $G$ is a connected graph with maximum degree $\Delta \geq 3$, then 
$$
Z(G) \le (\Delta - 2)\beta(G)+1,
$$
and this bound is sharp. 
\end{thma}
\proof Let $G$ be a connected graph with maximum degree $\Delta \geq 3$, let $C\subseteq V(G)$ be a minimum vertex cover of $G$, and let $X = V(G)\setminus C$. Note that $X$ is a maximum independent set, and so, every vertex in $C$ has at least one neighbor in $X$. Furthermore, since $G$ is connected, each vertex in $X$ has at least one neighbor in $C$. Color the vertices in $C$ blue and the vertices in $X$ white. Note that $C \cup X$ is a partition of $V(G)$, and so, if at any point $N_G[C]$ is colored blue, then all vertices in $G$ are colored blue (recall that every vertex in $X$ has a neighbor in $C$). Next, for each vertex $v \in C$ that has more than three neighbors in $X$, color blue a set of 
%\textcolor{red}{(since we are already talking about $v$ with more than 3 neighbors in $X$, shouldn't we just color $|N_G(v) \cap X|-3$ of its neighbors and be able to avoid using a maximum and potentially confusing reference to what is ``possible''?)} 
its neighbors in $X$ with cardinality $i$, where $i$ is the largest possible number that satisfies $i\leq d_G(v) - 3$. Let $S_0$ denote the resulting set of blue vertices and note that each vertex in $C$ has at most three white colored neighbors in $X$. Moreover, if $v\in S_0$ has three white colored neighbors, then all of the neighbors of $v$ are necessarily in $X$ (since we colored $d_G(v)-3$ neighbors of $v$ in $X$). Next let $\{C_1, C_2, C_3\}$ be a partition of $C$, where $C_i$ denotes the set of vertices with $i$ white colored neighbors in $X$. Since the set of vertices with three white neighbors in $X$ is an independent set in the subgraph induced by $C$, we note that $C_3$ is a set of isolates in the subgraph induced by $C$; that is, no vertex in $C_3$ is adjacent with another vertex of $C$. This fact is crucial to our argument, and we state this explicitly with the following observation.
\begin{ob}\label{ob:1}
The vertices in $C_3$ form a set of isolates in the subgraph induced by $C$, and this implies that no vertex in $C_3$ has a neighbor in $C$. 
\end{ob}

For the duration of this proof the sets $C_1$, $C_2$, and $C_3$ remain the same and do not alter at any step of our argument hereafter. With what follows, we construct a zero forcing set of $G$ by allowing the zero forcing process to start with the blue colored vertices in $S_0$, and then modify $S_0$ if all of $V(G)$ is not colored blue when the zero forcing process starting at $S_0$ terminates. Before moving forward with our construction, we first note the following inequality.
%\textcolor{red}{Since there may be vertices in $C$ with degree less than 3, this first sum may include negative terms.  However, the second sum is perfectly fine, so the right half of this first line just needs to be removed.}
\begin{equation}
\label{eq1}
\begin{array}{lcl}
 |S_0| &  \leq & \displaystyle{|C| + \sum_{v\in C}(\Delta-3)} \1 \\
& = & |C| + (\Delta -3) |C| \1 \\
& = & (\Delta-2)|C| \1 \\
& = & (\Delta-2)\beta(G).
\end{array}
\end{equation}

We now test if $S_0$ is a zero forcing set of $G$. Specifically, color all vertices in $S_0$ blue, color all other vertices white, and let the zero forcing process start and continue until no further color changes are possible. If all of $V(G)$ becomes colored blue, then $S_0$ is a zero forcing set of $G$. This observation together with inequality \eqref{eq1} imply
\[
Z(G) \le |S_0| \le (\Delta-2)\beta(G),
\]
which proves the proposed inequality. Thus, we will suppose $S_0$ is not a zero forcing set of $G$, since otherwise we have proven the desired inequality. 

Let $B_0$ be the set of blue colored vertices after termination of the zero forcing process starting at $S_0$ and let $W_0 = V(G)\backslash B_0$; in other words, $W_0$ is the set of white colored vertices at the termination of the zero forcing process starting with $S_0$. Since every vertex in $C_1$ initially has exactly one white colored neighbor, and this neighbor is in $X$, each vertex in $C_1$ can force its neighbor in $X$ to become colored blue at the first iteration of the zero forcing process. This implies that $N_G[C_1]$ is currently colored blue. Next observe that since $W_0$ is non-empty and since $N_G[C_1]$ is currently colored blue, there must be vertices in $C_2 \cup C_3$ that have (two or three) white colored neighbors in $W_0$. If there is a vertex $v_0\in C_2$ with two neighbors in $X\cap W_0$, say $w_0$ and $z_0$,  
%Let $C'_0$ be the set of vertices in $C$ that have all their neighbors in $B_0$. %Note that $w$ and $z$ cannot be adjacent to any vertices in $C'$. 
then we modify $S_0$ by replacing $v_0$ with $w_0$: 
$$
S_1 = (S_0\setminus \{v_0\})\cup \{w_0\}.
$$

We next test $S_1$ as a possible zero forcing set of $G$. Specifically, we color all vertices in $S_1$ blue, color all other vertices white, and let the zero forcing process start and continue until no further color changes are possible. Note that in $S_1$, the only initially white colored neighbor of $w_0$ is $v_0$, which implies $w_0$ can initially force $v_0$ to become colored blue. After $w_0$ forces $v_0$ to become colored blue, the only white colored neighbor of $v_0$ will be $z_0$. Thus, $v_0$ will then force $z_0$ to become colored blue. These color changes ensure that $N_G[v_0]$ will become colored blue by the zero forcing process starting at $S_1$. Moreover, after $N_G[v_0]$ becomes colored blue, any vertices in $C_2$ adjacent with $w_0$ or $z_0$ will have at most one white colored neighbor and will force this neighbor to become colored blue. This implies that eventually any vertex in $C_2$ adjacent with $w_0$ or $z_0$ will at some step of the zero forcing process have a blue colored closed neighborhood. After $N_G[v_0]$ becomes colored blue and the closed neighborhood of any vertex in $C_2$ adjacent with either $v_0$ or $z_0$ becomes colored blue, the set of currently blue colored vertices is a superset of $S_0$, and so, all of the color changes that occurred in $S_0$ may still occur in $S_1$, and in the same sequence. %Moreover, any vertices in $C_2$ adjacent to $w$ may also force, since any such vertex will have only one white neighbor in $X$. If all of $V(G)$ has become colored blue after the zero forcing process terminates, then $S$ is a zero forcing set of $G$. This observation together with inequality (1) imply 
If $S_1$ is a zero forcing set of $G$, then since $S_1$ has the same cardinality as $S_0$, inequality \eqref{eq1} implies that 
\[
Z(G) \le |S_1| \le (\Delta-2)\beta(G),
\]
which proves the proposed inequality. Hence, we will suppose $S_1$ is not a zero forcing set of $G$, since otherwise we have proven the desired inequality. 

Let $B_1$ be the set of blue colored vertices obtained after termination of the zero forcing process starting at $S_1$ and let $W_1 = V(G)\setminus B_1$. Then, and as before, there must be vertices in $C_2 \cup C_3$ that have (two or three) white colored neighbors in $X\cap W_1$. However, there is an additional condition that no such vertex in $C_2$ can be adjacent with either $w_0$ or $z_0$, since as described above, this vertex would have no white colored neighbors. Observing this, next suppose $v_1\in C_2$ is a vertex with two currently white colored neighbors in $X\cap W_1$, say $w_1$ and $z_1$, then we may modify $S_1$ by replacing $v_1$ with $w_1$: 
$$
S_2 = (S_1\setminus \{v_1\})\cup \{w_1\}.
$$
We next test $S_2$ as a possible zero forcing set of $G$. That is, we color all vertices in $S_2$ blue, color all other vertices white, and let the zero forcing process start and continue until no further color changes are possible. Since $N_G[v_0]$ was colored blue by starting the zero forcing process at $S_1$, $w_1$ is not adjacent to $v_0$. This implies that all of the neighbors of $w_1$ in $C \backslash \{v_1\}$ are initially colored blue. Thus, the only initially white colored neighbor of $w_1$ is $v_1$. Hence, $w_1$ initially may force $v_1$ to become colored blue. After $v_1$ becomes colored blue, we arrive at a blue colored superset of $S_1$, and so, all of the color changes that occurred in $S_1$ may still occur in $S_2$, and in the same sequence (even if $v_0$ is adjacent to $v_1$). After these color changes occur, the only possibly white colored neighbor of $v_1$ would be $z_1$. Thus, $v_1$ may force $z_1$ to become colored blue. This implies that $N_G[v_1]$ will also become colored blue by the zero forcing process starting at $S_2$. If $S_2$ is a zero forcing set of $G$, then since $S_2$ has the same cardinality as $S_0$, inequality \eqref{eq1} implies that 
\[
Z(G) \le |S_2| \le (\Delta-2)\beta(G),
\]
which is the proposed inequality. Hence, we will suppose $S_2$ is not a zero forcing set of $G$, since otherwise we have proven the desired inequality.

We next repeat this process until we obtain a set $S_k$ with cardinality equal to $S_0$, where after allowing the zero forcing process to start from $S_k$ and continue until no further color changes are possible, we are assured that the \textit{closed neighborhood} of vertices in $C_1 \cup C_2$ are eventually colored blue after the zero forcing process terminates. If $S_k$ is a zero forcing set of $G$, then since $S_k$ has the same cardinality as $S_0$, inequality \eqref{eq1} implies
\[
Z(G) \le |S_k| \le (\Delta-2)\beta(G),
\]
which proves the proposed inequality. Hence, we will suppose $S_k$ is not a zero forcing set of $G$, since otherwise we have proven the desired inequality.

Starting from the set $S_k$, let the zero forcing process start and continue until no further color changes are possible. Note that after termination, $N_G[C_1 \cup C_2]$ must be colored blue. Since we are assuming not all of $V(G)$ has become colored blue starting from the blue colored set $S_k$, it must be the case that there is some vertex in $C_3$ with two or three white colored neighbors in $X$ after the zero forcing process terminates. Next recall Observation~\ref{ob:1}, which states that $C_3$ is a set of isolates in the subgraph induced by $C$, and so, no vertex in $C_3$ is adjacent with any other vertex of $C$. Since $G$ is connected, this observation implies that either $C \cup X = C_3 \cup X$ is a bipartition of $V(G)$, or $C_1 \cup C_2 \neq \emptyset$. With what follows we consider both of these cases. 

We first suppose that $C_1 \cup C_2 \neq \emptyset$. This supposition together with $G$ being connected implies that there is at least one vertex in $C_3$ that shares a neighbor in $X$ with a vertex of $C_1 \cup C_2$ (since otherwise there would be no path from vertices in $C_3$ to vertices in $C_1\cup C_2$). Moreover, we know that $N_G[C_1 \cup C_2]$ has been colored blue. It follows that at least one vertex in $C_3$ has a neighbor that became colored blue by either a forcing vertex in $C_1\cup C_2$, or by the color modifications used to construct $S_k$. Furthermore, if any two neighbors of a vertex in $C_3$ have been colored blue, then this vertex may then force its white colored neighbor to become colored blue, which in-turn may force a neighbor of another vertex in $C_3$ to become colored blue. 

Since not all all of $V(G)$ has become colored blue, the above observations imply that there is at least one vertex $a_0\in C_3$ with exactly one neighbor $b_0 \in X$ that has been either forced to become colored blue by the zero forcing process, or has been colored blue by our construction of $S_k$. Let $c_0$ and $d_0$ denote the white colored neighbors of $a_0$ after termination of the zero forcing process starting at $S_k$. Since $c_0$ has yet to be colored blue, we know that $c_0$ is not adjacent with any vertices of $C_1\cup C_2$. This implies that $c_0$ is only adjacent with vertices of $C_3$, which in-turn implies that all of the neighbors of $c_0$ are initially colored blue in $S_k$. Moreover, no neighbor of $c_0$ could have been a forcing vertex during the zero forcing process since otherwise $c_0$ would be colored blue. We next modify $S_k$ by replacing $a_0$ with $c_0$:
\[
S_k^1 = (S_k \setminus \{a_0\}) \cup \{c_0\}
\]
We next test $S_k^1$ as a possible zero forcing set of $G$. That is, we color all vertices in $S_k^1$ blue, color all other vertices white, and let the zero forcing process start and continue until no further color changes are possible. Note that $c_0$ could not have been adjacent with a vertex of $C_1 \cup C_2$, since otherwise $c_0$ would have been colored blue by $S_k$. Thus, the only possible initially white colored neighbor of $c_0$ is $a_0$. It follows that $c_0$ may initially force $a_0$ to become colored blue. After $c_0$ forces $a_0$ to become colored blue, we arrive at a blue colored superset of $S_k$, and so, all color changes which occurred during the zero forcing process starting from the blue set of vertices in $S_k$ will now happen in $S_k^1$. This implies that eventually the neighbor $b_0$ of $a_0$ will become colored blue. Whenever $b_0$ is colored blue (either by initially being blue colored or forced to become colored blue), the only white colored neighbor of $a_0$ will be $d_0$. Thus, $a_0$ may force $d_0$ to become colored blue. Hence, we are assured that the closed neighborhood of $a_0$ has become colored blue by the time the zero forcing process starting at $S_k^1$ terminates. If $S_k^1$ is a zero forcing set of $G$, then since $S_k^1$ has the same cardinality as $S_0$, inequality \eqref{eq1} implies
\[
Z(G) \le |S_k^1| \le (\Delta-2)\beta(G),
\]
which proves the proposed inequality. Hence, we will suppose $S_k^1$ is not a zero forcing set of $G$, since otherwise we have proven the desired inequality. 

Let $B_k^1$ be the set of blue colored vertices obtained after termination of the zero forcing process starting at $S_k^1$ and let $W_k^1 = V(G)\setminus B_k^1$. Since not all of $V(G)$ has become colored blue, the aforementioned observations on vertices in $C_3$ imply that there is at least one vertex $a_1\in C_3$ with exactly two white colored neighbors in $X\cap W_k^1$ and exactly one neighbor, say $b_1 \in X\cap B_k^1$, such that $b_1$ has been forced to become colored blue by the zero forcing process starting at $S_k^1$. Let $c_1$ and $d_1$ be the white colored neighbors of $a_1$ in $X\cap W_k^1$. Since $c_1$ has yet to be colored blue, we know that $c_1$ is not adjacent with any vertices of $C_1\cup C_2$. Next observe that neither $c_1$ nor $d_1$ could be adjacent with a vertex that has previously performed a force, since otherwise they would have been forced to become colored blue and they are both colored white. With this observation we next modify $S_k^1$ by replacing $a_1$ with $c_1$:
\[
S_k^2 = (S_k^1 \setminus \{a_1\}) \cup \{c_1\}
\]

We next test $S_k^2$ as a possible zero forcing set of $G$. That is, we color all vertices in $S_k^2$ blue, color all other vertices white, and let the zero forcing process start and continue until no further color changes are possible. Note that $c_1$ could not have been adjacent with a vertex of $C_1 \cup C_2$, nor be adjacent with $a_0$ since otherwise $c_1$ would have been colored blue by the zero forcing process starting at $S_k^1$. Thus, the only possible initially white colored neighbor of $c_1$ is $a_1$. It follows that $c_1$ may initially force $a_1$ to become colored blue. After $c_1$ forces $a_1$ to become colored blue, we arrive at a blue colored superset of $S_k^1$, and so, all color changes which occurred during the zero forcing process starting from the blue set of vertices in $S_k^1$ will now happen in $S_k^2$. This implies that eventually the neighbor $b_1$ of $a_1$ will become colored blue by starting the zero forcing process at $S_k^2$. After $b_1$ is colored blue (either by being forced to become colored blue, or by initially being colored blue) the only white colored neighbor of $a_1$ will be $d_1$. Thus, $a_1$ may force $d_1$ to become colored blue. Hence, we are assured that the closed neighborhood of $a_1$ is colored blue after the zero forcing process terminates starting at $S_k^2$. If $S_k^2$ is a zero forcing set of $G$, then since $S_k^2$ has the same cardinality as $S_0$, inequality \eqref{eq1} implies
\[
Z(G) \le |S_k^2| \le (\Delta-2)\beta(G),
\]
which proves the proposed inequality. Hence, we will suppose $S_k^2$ is not a zero forcing set of $G$, since otherwise we have proven the desired inequality.

We next repeat this process until we obtain a set $S_k^j$ with cardinality equal to $S_0$, where after allowing the zero forcing process to start from $S_k^j$ and continue until no further color changes are possible, we are assured that the closed neighborhoods of vertices in $C_1 \cup C_2 \cup C_3$ are colored blue by the zero forcing process starting at $S_k^j$.; that is, $N_G[C_1 \cup C_2 \cup C_3]$ is colored blue. This implies that $S_k^j$ is a zero forcing set of $G$. Since $S_k^j$ has the same cardinality as $S_0$, inequality \eqref{eq1} implies
\[
Z(G) \le |S_k^j| \le (\Delta-2)\beta(G).
\]

At this point of our argument we have proven that if $C_1 \cup C_2 \neq \emptyset$, then 
$$
Z(G) \le (\Delta-2)\beta(G) < (\Delta-2)\beta(G) + 1.
$$
Thus, we will assume $C_1 \cup C_2 = \emptyset$, since otherwise we have proven the desired inequality. Assuming $C_1 \cup C_2 = \emptyset$ and recalling Observation~\ref{ob:1}, it must be the case that $C \cup X = C_3 \cup X$ is a bipartition of $V(G)$. Moreover, under this supposition, the fact that every vertex of $C_3$ has exactly three white neighbors implies that no color changes may  occur starting from $S_k$, and so, $S_k = S_0$. Thus, in our process of constructing a zero forcing set of $G$, we are starting from the set $S_0$. To circumvent the problem of no vertex being able to initially force, we select an arbitrary vertex in $C_3$, say $a_0$, with neighbors $b_0$, $c_0$, and $d_0$, in $X$. Next we \emph{greedily} add the vertex $b_0$ to our initial set $S_0$. This modification of $S_0$ is given by:
\[
S_1 = S_0 \cup \{b_0\}
\]
Note that we have added exactly one vertex to $S_0$. Thus, inequality \eqref{eq1} implies
\begin{equation}\label{eq2}
|S_1| = |S_0| + 1 \le (\Delta-2)\beta(G) + 1.
\end{equation}

Next color all vertices in $S_1$ blue and color all other vertices white. Starting from the set $S_1$, observe that no color changes may occur during the zero forcing process (since no blue colored vertex has exactly one white colored neighbor). However, we also observe that $a_0 \in C_3$ is colored blue and its neighbor $b_0\in X$ is also colored blue. Next recall that $c_0$ and $d_0$ are the white colored neighbors of $a_0$, and then modify $S_1$ by replacing $a_0$ with $c_0$:
\[
S_1^1 = (S_1 \setminus \{a_0\}) \cup \{c_0\}
\]

We next test $S_1^1$ as a possible zero forcing set of $G$. That is, we color all vertices in $S_1^1$ blue, color all other vertices white, and let the zero forcing process start and continue until no further color changes are possible. Note that $c_0$ is colored blue and all neighbors of $c_0$, other than $a_0$, are also initially colored blue. It follows that $c_0$ may initially force $a_0$ to become colored blue. After $c_0$ forces $a_0$ to become colored blue, we next observe that $d_0$ is the only white colored neighbor of $a_0$. This implies that on the second iteration of the zero forcing process starting at $S_1^1$, the vertex $a_0$ may force $d_0$ to become colored blue. Hence, we are assured that the closed neighborhood of $a_0$ has become colored blue after the zero forcing process terminates starting at $S_1^1$. If $S_1^1$ is a zero forcing set of $G$, then since $S_1^1$ has the same cardinality as $S_1$, inequality~\eqref{eq2} implies
\[
Z(G) \le |S_1^1| \le (\Delta-2)\beta(G) + 1,
\]
which proves the proposed inequality. Hence, we will suppose $S_1^1$ is not a zero forcing set of $G$, since otherwise we have proven the desired inequality.

Starting from the set $S_1^1$, let the zero forcing process start and continue until no further color changes are possible. Note that not all of $V(G)$ has become colored blue and the closed neighborhood of $a_0$ is colored blue. Since $G$ is connected and Observation~\ref{ob:1} implies that $C_3$ is an independent set, it must be the case that a neighbor of $a_0$ in $X$ is adjacent with a vertex different than $a_0$ in $C_3$ (otherwise there would not be a path from $a_0$ to other vertices in $C_3$). Since the closed neighborhood of $a_0$ has become colored blue, we are assured that at least one vertex in $C_3$, different than $a_0$, has had at least one neighbor in $X$ become colored blue (either by being forced to become colored blue or by our construction of $S_1^1)$. Moreover, if any vertex in $C_3$ has exactly two blue colored neighbors, then this vertex could then force its one white colored neighbor to become colored blue. Since not all of $V(G)$ has become colored blue, these observations imply that there is at least one vertex in $C_3$ which has exactly two white colored neighbors in $X$ and exactly one blue colored neighbor in $X$ after the zero forcing process starting from $S_1^1$ terminates. Let $a_1 \in C_3$ be one such vertex, let $b_1 \in X$ be the one blue colored neighbor of $a_1$ in $X$, and let $c_1$ and $d_1$ be the white colored neighbors of $a_1$ in $X$. We next modify $S_1^1$ by replacing $a_1$ with $c_1$:
\[
S_1^2 = (S_1^1 \setminus \{a_1\}) \cup \{c_1\}
\]

We next test $S_1^2$ as a possible zero forcing set of $G$. That is, we color all vertices in $S_1^2$ blue, color all other vertices white, and let the zero forcing process start and continue until no further color changes are possible. Note that $c_1$ could not have been adjacent with $a_0$ since otherwise $c_1$ would have been colored blue by the zero forcing process starting from $S_1^1$. Thus, the only possible initially white colored neighbor of $c_1$ is $a_1$. It follows that $c_1$ may initially force $a_1$ to become colored blue. After $c_1$ forces $a_1$ to become colored blue, we arrive at a blue colored superset of $S_1^1$, and so, all color changes which occurred during the zero forcing process starting from $S_1^1$ may now happen in $S_1^2$. This implies that eventually $b_1$ will become colored blue. After $b_1$ is colored blue, the only white colored neighbor of $a_1$ will be $d_1$. Thus, $a_1$ may then force $d_1$ to become colored blue. Hence, we are assured that the closed neighborhood of $a_1$ has become colored blue. If $S_1^2$ is a zero forcing set of $G$, and since $S_1^2$ has the same cardinality as $S_1$, inequality (2) implies
\[
Z(G) \le |S_1^2| \le (\Delta-2)\beta(G) + 1,
\]
which proves the proposed inequality. Hence, we will suppose $S_1^1$ is not a zero forcing set of $G$, since otherwise we have proven the desired inequality.

We next repeat this process until we obtain a set $S_1^j$ with cardinality equal to $S_1$, where after allowing the zero forcing process to start from $S_1^j$ and continue until no further color changes are possible, we are assured that the closed neighborhood of vertices in $C = C_3$ are colored blue after the zero forcing process terminates. It follows that $S_1^j$ is a zero forcing set of $G$, then since $S_1^j$ has the same cardinality as $S_1$, inequality~\eqref{eq2} implies
\[
Z(G) \le |S_2^j| \le (\Delta-2)\beta(G) + 1,
\]
which proves the proposed inequality. In particular, we have now proven that if $C_3 \cup X$ is a bipartition of $V(G)$, then 
\[
Z(G) \le (\Delta-2)\beta(G) + 1.
\]

Since either $C_1 \cup C_2 \neq \emptyset$, or $C_3 \cup X$ is a bipartition of $V(G)$, we have now shown that if $G$ is a connected graph with maximum degree $\Delta \ge 3$, then $Z(G) \le (\Delta-2)\beta(G) + 1$. To see this bound is sharp see Section~\ref{sec:sharp2}. \qed

\subsection{Sharp Examples of Theorem \ref{theorem2}}\label{sec:sharp2}
\begin{figure}[htb]
\begin{center}
\includegraphics{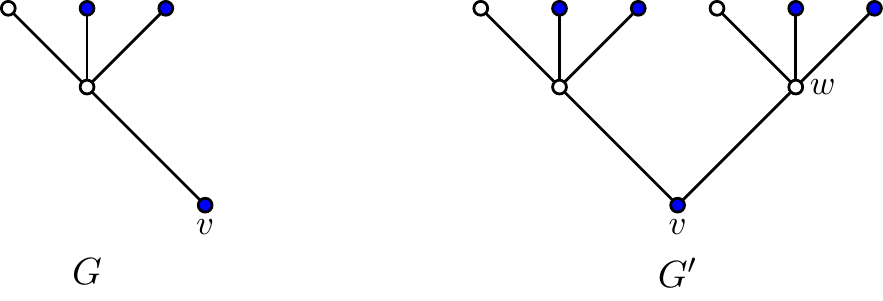}
\end{center}
\caption{The star graph $G = K_{1, 4}$ and the graph $G'$ obtained by performing a 3-LSVA on $K_{1, 4}$. Minimum zero forcing sets of $G$ and $G'$ shown in blue.}
\label{fig:lsva}
\end{figure}

In this section we construct a broad infinite family of graphs satisfying Theorem \ref{theorem2} with equality. We begin by noting that for the star graph $K_{1, n-1}$, $n\ge 4$, it holds that $Z(K_{1, n-1})= n-2 =(\Delta(K_{1, n-1})-2)\beta(K_{1, n-1}) + 1$. 
We next define a general operation called \emph{$k$-leaf support vertex addition} on a graph $G$, abbreviated $k$\emph{-LSVA}. Given a graph $G$ with maximum degree $\Delta$, we define $k$-LSVA on $G$ to be the process of attaching to a vertex $v\in V(G)$ with degree $d_G(v) \le k -1$ a new vertex $w$, and then attaching $k$ leaves to $w$; see Figure~\ref{fig:lsva} for an illustration. 

Using $k$-LSVA we may iteratively generate an infinite family of graphs with maximum degree $\Delta \ge 3$ and arbitrary diameter which satisfy $Z(G) = (\Delta-2)\beta(G) + 1$. Before giving this construction, we first prove the following useful lemma.

\begin{lem}\label{prop:tech}
If $G$ is a graph with $\Delta(G) \ge 3$ satisfying $Z(G) = (\Delta(G)-2)\beta(G) + 1$ and if $G'$ is any graph obtained by applying $(\Delta(G) - 1)$-LSVA on $G$, then $Z(G') = (\Delta(G')-2)\beta(G') + 1$. 
\end{lem}

\proof Let $G$ be a graph with maximum degree $\Delta \ge 3$ satisfying $Z(G) = (\Delta-2)\beta(G) + 1$, and let $v \in V(G)$ be an arbitrary vertex with degree $d_G(v) \le \Delta - 1$. Next perform $(\Delta - 1)$-LSVA on $G$, by attaching a vertex $w$ to $v$ and then attaching $\Delta - 1$ leaves to $w$; denote the resulting graph by $G'$. Next observe that since $\Delta \ge 3$, the leaves of $w$ are contained in every maximum independent set in $G'$. This implies that $w$ is contained in every minimum vertex cover of $G'$. Thus, $\beta(G') = \beta(G) + 1$. Next note that if any minimum zero forcing set of $G$ together with $\Delta-2$ leaf neighbors of $w$ is a zero forcing set of $G'$, then this set is a minimum zero forcing set of $G'$.  Now let $S\subseteq V(G)$ be a fixed minimum zero forcing set of $G$, let $W = \{w_1, \dots, w_{\Delta-2}\}$ be $\Delta-2$ arbitrary leaf neighbors of $w$, and let $S' = S \cup W$. Color the vertices of $S'$ blue and all other vertices color white. With this coloring observe that one of the leaf neighbors of $w$, say $w_1$, may initially force $w$ to become colored blue. Next observe that once $w$ is colored blue every vertex in $S$ is blue and no vertex in $S$ has a white neighbor outside of $G$, and so, every color change occurring in $G$ starting from $S$ will also occur. Since $S$ was a zero forcing set of $G$, it follows that $v$ will eventually become colored blue. Once $v$ is colored blue, $w$ will have exactly one white colored neighbor (a leaf). It follows that $w$ may then force its white colored leaf neighbor to become colored blue. Since $S$ is assured to color all vertices in $V(G)$ blue and since $N_{G'}[w]$ has become colored blue it must be the case that $S'$ is a zero forcing set of $G'$. Moreover, $S'$ is a minimum zero forcing set of $G'$ as described earlier. Thus, $Z(G') = Z(G) + (\Delta - 2)$. Hence, 
\begin{equation*}
\begin{array}{lcl}
 Z(G') &  = & Z(G) + (\Delta - 2) \1 \\
& = & (\Delta - 2)\beta(G) + 1 + (\Delta - 2) \1 \\
& = & (\Delta - 2)[\beta(G') - 1] + \Delta - 1 \1 \\
& = & \Delta\beta(G') - \Delta - 2\beta(G') + 2 + \Delta - 1 \1 \\
& = & (\Delta-2)\beta(G') + 1.
\end{array}
\end{equation*}
\qed 

\begin{figure}[htb]
\begin{center}
\includegraphics{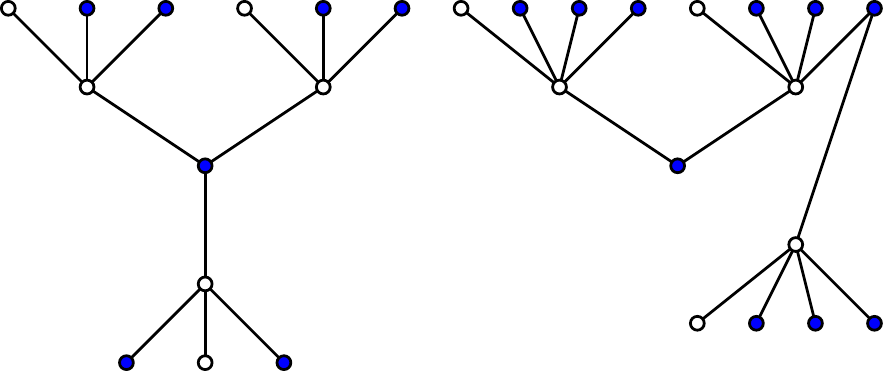}
\end{center}
\caption{Trees $G \in \mathcal{T}^*$ satisfying $Z(G) = (\Delta - 2)\beta(G) + 1$. Minimum zero forcing sets shown in blue.}
\label{fig:general}
\end{figure}

Using Lemma~\ref{prop:tech}, we  construct an infinite family of trees with arbitrary maximum degree $\Delta \ge 3$ and arbitrarily large diameter. Let $k\ge 3$ and let $\mathcal{T}^*$ be the set of graphs obtained by starting with $K_{1, k}$ and applying as many $(k-1)$-LSVA's as wanted. Since $K_{1, k}$ satisfies Theorem \ref{theorem2} with equality, Lemma~\ref{prop:tech} implies that all graphs in $\mathcal{T}^*$ also satisfy Theorem \ref{theorem2} with equality. See Figure~\ref{fig:general} for illustrations of such graphs.

Finally, we consider the special case of Theorem \ref{theorem2} when $\Delta(G)=3$. In particular, this yields the following elegant corollary.
\begin{cor}\label{cor:main}
If $G$ is a connected graph with maximum degree $\Delta = 3$, then
\[
Z(G) \le \beta(G) + 1,
\]
and this bound is sharp. 
\end{cor}

\begin{figure}[htb]
\begin{center}
\includegraphics{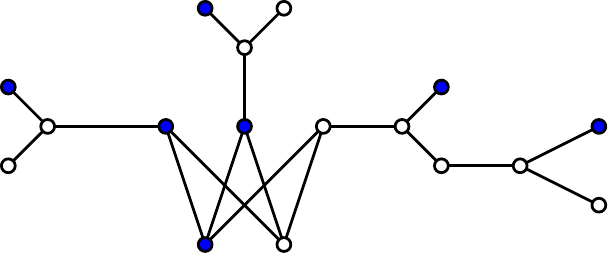}
\end{center}
\caption{A graph $G \in \mathcal{G}^*$ with maximum degree $\Delta = 3$ which satisfies $Z(G) = \beta(G) + 1$; a minimum zero forcing set shown in blue.}
\label{fig:general2}
\end{figure}

\begin{figure}[htb]
\begin{center}
\includegraphics{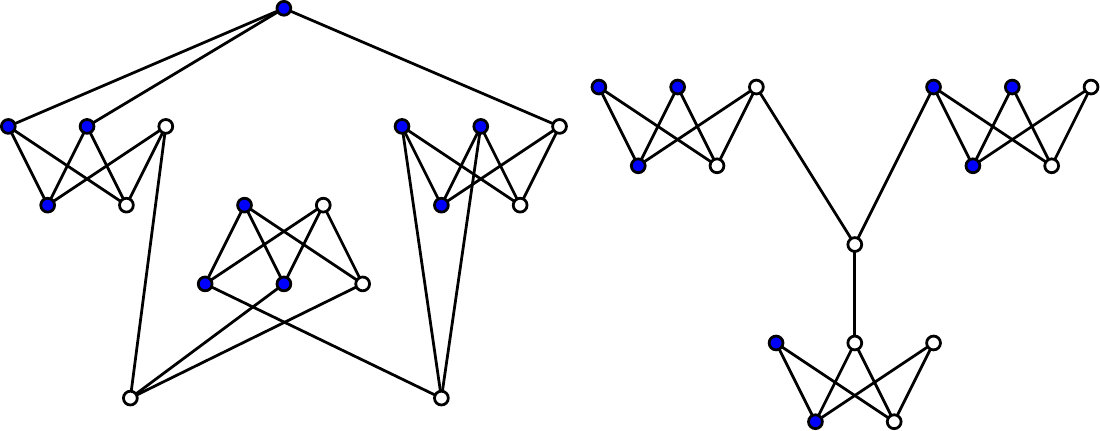}
\end{center}
\caption{Graphs with maximum degree $\Delta = 3$ that satisfy $Z(G) = \beta(G) + 1$. Minimum zero forcing sets shown in blue.}
\label{fig:general4}
\end{figure}

We next present an infinite family of graphs that satisfies Corollary~\ref{cor:main} with equality. Let $\mathcal{G}^*$ be the set of graphs obtained by starting with $K_{2, 3}$ and applying an arbitrary number of $2$-LSVA's. Since $K_{2, 3}$ satisfies Theorem \ref{theorem2} with equality, Lemma~\ref{prop:tech} implies that all graphs in $\mathcal{G}^*$ also satisfy Corollary~\ref{cor:main} with equality. See Figure~\ref{fig:general2} for an illustration of one such graph. Note that the graphs in $\mathcal{G}^*$ are not trees; in contrast, for values of $\Delta$ greater than three, we have not found any graphs that are not trees that satisfy Theorem~\ref{theorem2} with equality. The graphs in $\mathcal{G}^*$ and $\mathcal{T}^*$ are not the only graphs that satisfy Theorem~\ref{theorem2} (or Corollary~\ref{cor:main}) with equality. We have found other more elaborate constructions of such graphs (see, e.g., Figure~\ref{fig:general4}) but we do not have a complete characterization of their structure.

\section{Conclusion}
\label{conclusion}
In this paper we have shown that machine learning and automated conjecturing can produce new and interesting theorems that in turn motivate other significant results. The initial conjecture of \emph{TxGraffiti} that we considered also led to the discovery of a new polynomial time algorithm to obtain a zero forcing set in claw free graphs, and also to constructions of infinite families of graphs in which the zero forcing number equals the vertex cover number. 

One direction for future work is to obtain more precise characterizations of the graphs in which the obtained bounds hold with equality. In particular, we believe that the family of trees $\mathcal{T}^*$ constructed in Section~\ref{sec:sharp2} are the only trees satisfying Theorem \ref{theorem2} with equality. We pose this formally with the following question. 
\begin{quest}
Let $T$ be a tree with maximum degree $\Delta \ge 3$. Is it true that 
$$
Z(T) = (\Delta - 2)\beta(T) + 1,
$$ 
if and only if $T\in \mathcal{T}^*$? 
\end{quest}

More generally, we also suggest the following question.
\begin{quest}
Which connected graphs $G$ with $\Delta(G) \ge 3$ satisfy $Z(G) = (\Delta - 2)\beta(G) + 1$?
\end{quest}

Next recall that the vertex cover number of $G$ and the independence number of $G$ are related by the equation $\beta(G) = n(G) - \alpha(G)$. This relationship and the results presented in this paper encourage further study of the following open problem, which was also conjectured by \emph{TxGraffiti}.
\begin{conj}[TxGraffiti 2017]\label{conj:alpha}
If $G \ne K_4$ is a connected graph with $\Delta(G) \le 3$, then
$$
Z(G) \le \alpha(G) + 1,
$$ 
and this bound is sharp.
\end{conj}

When considering Conjecture~\ref{conj:alpha}, we note that a partial result of its statement has been given by Davila and Henning in~\cite{DaHe19c}; namely, if $G\neq K_4$ is a claw-free and 3-regular graph, then $Z(G) \le \alpha(G) + 1$. Furthermore, the proof of this result is similar to the iterative nature of the proofs of Theorem \ref{thm1} and Theorem \ref{theorem2} given in this paper. Thus, we suspect that some combination of these techniques and ideas may eventually lead to a proof of Conjecture~\ref{conj:alpha}.

\end{document}